\documentclass[11pt]{article}
\usepackage{amsfonts}

\newcommand{\sect}[1]{\setcounter{equation}{0}\section{#1}}
\newcommand{\subsect}[1]{\subsection{#1}}

\def\be{\begin{equation}}
\def\ee{\end{equation}}
\def\bea{\begin{eqnarray}}
\def\eea{\end{eqnarray}}

\def\1{\'{\i}}                           
\def\R{{\mathbb R}}

\def\bdiez{\alpha}
\def\btres{\xi}
\def\bcuatro{\nu}
\def\bcinco{\beta_1}
\def\bseis{\beta_2}
\def\bdoce{\beta_3}
\def\bocho{\beta_4}
\def\bonce{\beta_5}
\def\buno{\beta_6}

\def\auno{\btres}
\def\ados{\bdoce}
\def\asiete{\bseis}
\def\adiez{\bcinco}

\def\aauno{a_1}
\def\aados{a_2}
\def\aatres{a_3}
\def\aaseis{a_4}
\def\aasiete{a_5}
\def\aadiez{a_{6}}

\def\gal{\overline{{\cal G}}}

\def\masa{{\cal M}}

\begin{document}

\thispagestyle{empty}

 \
\hfill\ J.Phys. A 33 (2000) 3431-3444

\ 
\vspace{2cm}

\begin{center}
{\LARGE{\bf{Quantum (1+1) extended Galilei algebras:}}}

{\LARGE{\bf{from Lie bialgebras to quantum $R$-matrices}}}

{\LARGE{\bf{and integrable systems}}} 
\end{center}

\bigskip\bigskip

\begin{center} Angel Ballesteros$^\dagger$,  Enrico
Celeghini$^\ddagger$ and Francisco J. Herranz$^\dagger$ 
\end{center}

\begin{center} {\it {$^\dagger$ Departamento  de F\1sica,
Universidad de Burgos} \\   Pza. Misael Ba\~nuelos, 
E-09001 Burgos, Spain}
\end{center}

\begin{center} 
{\it {$^\ddagger$  Dipartimento di Fisica and  Sezione INFN,
Universit\'a di   Firenze\\ I-50125 Firenze, Italy}}
\end{center}

\bigskip\bigskip

\begin{abstract}
The Lie bialgebras of the (1+1) extended  Galilei algebra are
obtained and classified into four multiparametric families.
Their quantum deformations are obtained, together with the
corresponding deformed Casimir operators. For the coboundary
cases quantum universal $R$-matrices are also given.
Applications of the quantum extended Galilei algebras to
classical integrable systems  are  explicitly developed.
\end{abstract}

\newpage


\sect{Introduction}

The study of  Lie bialgebra structures  
 provides a primary classification of 
the zoo of possible quantum deformations of  a given Lie
algebra \cite{CP}. For simple Lie algebras, this
problem has been studied in \cite{BelDr,Stolin}; in
this case, all Lie bialgebras are of the coboundary type and
their classification reduces to obtain all 
constant solutions of the classical Yang--Baxter equation. 
During the last years, 
the classification of the Lie bialgebras (and, sometimes, of
the corresponding Poisson--Lie structures)  
for some non-simple Lie algebras with  physical
interest have been found. The results cover mainly low
dimensional cases: the Heisenberg--Weyl
$h_3$ or (1+1) Galilei algebra 
\cite{kuper,vero,heis,galileo}, the two-dimensional Euclidean
algebra \cite{euclideo}, the harmonic oscillator
$h_4$ algebra \cite{verob,osc}, the (1+1) extended Galilei
algebra \cite{anna} and the $gl(2)$ algebra \cite{kuperb,gl}.
For higher dimensions, only the (3+1) Poincar\'e
algebra was treated in \cite{zakr,woro}.

In this paper we classify the (1+1)  extended Galilei 
$\gal$ Lie bialgebras in order  to obtain the quantum
deformations associated to $\gal$, and to show how  these
deformed Hopf structures can be directly used in some
applications such as   integrable models and deformed
heat-Schr\"odinger equations. With this in mind,   in the
next section all the $\gal$ Lie bialgebras are casted into
four  multiparametric families which naturally follow by
considering if the central generator is either a primitive or
a non-primitive generator.  The coboundary cases are
identified and it is shown that they belong to the first
family of bialgebras.  We stress that a classification of the
inequivalent   $\gal$ Lie bialgebras  together with
the Poisson--Lie structures has been  obtained  by Opanowicz
\cite{anna} while their corresponding quantum deformations
have been constructed in \cite{annab}. However our 
classification  in   multiparametric families is well adapted
and more manageable in order to  construct  systematically 
the  quantum  $\gal$ algebras; this is performed in section 3
by applying the   formalism introduced by Lyakhovsky  and   
Mudrov \cite{LM,Lyak}. In particular, for each
multiparametric quantum $\gal$ algebra, we obtain the 
coproduct,  the  compatible commutation rules and the 
Casimirs.  Furthermore, both standard and non-standard
quantum universal
$R$-matrices are deduced  for the coboundary quantum $\gal$
algebras in  section 4. As an application, we  show in section
5 the classical completely integrable systems that can be
constructed from these quantum algebras. We end the paper with
some comments concerning a space discretization of the
heat-Schr\"odinger equation  with quantum Galilei symmetry.


\sect{Extended Galilei  bialgebras}

The (1+1) extended Galilei algebra $\gal$ is a
 four-dimensional real Lie algebra  generated by $K$ (boost),
$H$ (time translation), $P$ (space translation) and $M$ (mass
of a particle in a free kinematics). The   Lie brackets and
Casimir operators of $\gal$ are given by
\be
[K,H]=P \qquad [K,P]=M \qquad [H,P]=0 \qquad [M,\cdot\,]=0.
\label{bc}
\ee
\be
{\cal C}_1=M\qquad {\cal C}_2=P^2 - 2 M H.
\label{bd}
\ee

In order to  obtain the Lie bialgebras associated to $\gal$
we have to find the most general  cocommutator $\delta:\gal
\rightarrow \gal 
\otimes\gal$ such that

\noindent i) $\delta$ is a 1-cocycle, i.e.,
\be
\delta([X,Y])=[\delta(X),\, 1\otimes Y+ Y\otimes 1] + 
[1\otimes X+ X\otimes 1,\, \delta(Y)]  \qquad \forall X,Y\in
\gal. 
\label{ba}
\ee
\noindent ii) The dual map $\delta^\ast:\gal^\ast\otimes
\gal^\ast \to \gal^\ast$ is a Lie bracket on $\gal^\ast$.

To begin with we consider a generic  linear combination (with
real coefficients) of skewsymmetric products of  the
generators $X_l$ of $\gal$:
\be
\delta(X_i)=f_i^{jk}\,X_j\wedge X_k  .
\label{be}
\ee
 By imposing
the cocycle condition (\ref{ba}) onto 
(\ref{be})  we find  the following  (pre)cocom\-mutator
which depends on nine parameters
$\{\bdiez,\btres,\bcuatro,\bcinco,
\bseis,\bdoce,\bocho,\bonce,\buno\}$:
\bea
&&\delta(K)=\buno K\wedge P+ \btres K\wedge M + \bcuatro
P\wedge H +
\bcinco P\wedge M + \bseis H\wedge M\cr
&&\delta(H)=\bonce K\wedge M -(\buno+\bdiez)  P\wedge H +
\bdoce P\wedge M + (\bocho - \btres) H\wedge M\cr
&&\delta(P)=\bocho P\wedge M + (\buno + \bdiez) H\wedge M\cr
&&\delta(M)=\bdiez P\wedge M .
\label{bf}
\eea
The Jacobi identities have to be imposed onto the dual map
$\delta^\ast$  in order to guarantee that this map defines 
Lie brackets. Thus we obtain the following set of equations:
\bea
&&\bdiez\bonce=0\qquad
\buno(\buno +\bdiez)=0\qquad
\bocho(\buno +\bdiez)=0\cr
&&\bcuatro(\btres -\bocho)=0\qquad
\bdiez(\btres -\bocho) - \bcuatro\bonce=0 .
\label{bg}
\eea
We solve the equations according  to the value of the
parameter $\bdiez$ since it characterizes the bialgebras with
primitive and non-primitive  mass $(\bdiez=0$ and $\bdiez\ne
0$, respectively).  In this way, it can be checked that the
general solution of (\ref{bg}) can be  splitted  into four
disjoint classes:

\noindent
{\em Family I: $M$ is a primitive generator}.

Ia)  $\bdiez=0$, $\buno=0$, $\bcuatro=0$ and  
$\{\btres,\bcinco,\bseis,\bdoce,\bocho,\bonce\}$ are arbitrary.

Ib) $\bdiez=0$, $\buno=0$, $\bcuatro\ne 0$, $\bocho=\btres$,
$\bonce=0$ and  $\{\btres,\bcinco,\bseis,\bdoce\}$ are
arbitrary.

\noindent
{\em Family II: $M$ is a non-primitive generator}.

IIa) $\bdiez\ne 0$, $\bonce=0$, $\buno=0$, $\btres=0$,
$\bocho=0$   and  
$\{\bcuatro,\bcinco,\bseis,\bdoce\}$ are arbitrary.

IIb) $\bdiez\ne 0$,  $\bonce=0$, $\buno= -\bdiez$,
$\bocho=\btres$ and 
$\{\btres,\bcuatro,\bcinco,\bseis,\bdoce\}$ are arbitrary.

 We recall that two Lie bialgebras $(\gal,\delta)$  and
$(\gal,\delta')$ are said to be equivalent if there exists
an automorphism ${\cal O}$ of $\gal$ such that
 $\delta'=({\cal O}\otimes {\cal O})\circ\delta\circ {\cal
O}^{-1}$.  The general  automorphism which preserves the 
commutation rules of $\gal$ (\ref{bc}) turns out to be
\be
\begin{array}{l}
K'=K+\lambda_1 H + \lambda_2 P + \lambda_3 M\cr
H'=H + \lambda_4 P + \lambda_5 M\cr
P'=P+\lambda_4 M\cr
M'=M 
\end{array}
\label{bbh}
\ee
where $\lambda_i$ are arbitrary real parameters. In what
follows we show how  this map enables us to simplify the
families of bialgebras  with some parameter different from
zero (Ib, IIa and IIb) by  removing  superflous parameters.

\noindent
$\bullet$ {Family Ib.} If we define 
\bea
&&K'=K\qquad 
H'=H-\frac {\bseis}{\bcuatro} P +\left(\frac
{\bcinco}{\bcuatro}+\frac {\bseis^2}{\bcuatro^2}\right) M\cr
&&P'=P-\frac {\bseis}{\bcuatro}M\qquad M'=M\qquad
\bdoce'=\bdoce - \frac{\bseis\btres}{\bcuatro}\qquad 
\bcuatro\ne 0 
\label{bh}
\eea
we obtain that the  cocommutators for the new generators
$X'$ are given by
\be
\begin{array}{l}
\delta(K')=\btres K'\wedge M' + \bcuatro P'\wedge H'
\qquad \delta(H')=\bdoce' P'\wedge M'\cr
\delta(P')=\btres P'\wedge M'\qquad \delta(M')=0 .
\end{array}
\ee
Therefore the parameters $\bcinco$ and $\bseis$ have been
removed from the cocommutators, so that   this family depends
on three parameters  $\{\bcuatro\ne 0,\btres,\bdoce\}$ with
$\bocho=\btres$.

\noindent
$\bullet$ {Family IIa.} We consider the
automorphism defined by:
\bea
&&K'=K+ \frac {\bcuatro}{\bdiez} H -\frac {\bseis}{\bdiez} P
-\left(\frac{\bcinco}{\bdiez}
+\frac{\bcuatro\bdoce}{\bdiez^2}\right)M
\cr &&H'= H -\frac{\bdoce}{2\bdiez} M\qquad P'=P\qquad
M'=M\qquad \bdiez\ne 0.
\label{bi}
\eea
The  cocommutators reduce to
\be
\delta(K')=0\qquad  \delta(H')=-\bdiez P'\wedge H'\qquad
 \delta(P')=\bdiez
H'\wedge M'\qquad
\delta(M')=\bdiez P'\wedge M'  .
\ee
Hence the parameters
$\{\bcuatro,\bcinco,\bseis,\bdoce\}$ have been reabsorbed and
this family depends on a single parameter $\bdiez\ne 0$.

\noindent
$\bullet$ {Family IIb.} In this  case there are three
superflous parameters $\{\bcuatro,\btres,\bdoce\}$ which
disappear when we define
\bea
&& K'=K+\frac{\bcuatro}{\bdiez} H\qquad 
H'=H-  \frac{\btres}{\bdiez} P 
+\left( \frac{\btres^2}{\bdiez^2}- 
\frac{\bdoce}{\bdiez} \right)
M\cr
&& P'=P -  \frac{\btres}{\bdiez} M\qquad M'=M 
\qquad
\bdiez \ne 0 \cr
&&\bcinco'=\bcinco +
 \frac {\bseis \btres}{\bdiez} +
 \frac {\bdoce \bcuatro}{\bdiez} - 
 \frac {\bcuatro \btres^2}{\bdiez^2}
\qquad
\bseis'=\bseis - \frac {\bcuatro \btres}{\bdiez}.
\label{bk}
\eea
The resulting cocommutators read
\be
\begin{array}{l}
\delta(K')=-\bdiez K'\wedge P'   + \bcinco' P'\wedge M'
+\bseis' H'\wedge M' \cr
 \delta(H')=0\qquad \delta(P')=0\qquad
\delta(M')=\bdiez P'\wedge M'  
\end{array}
\ee
and this equivalence of bialgebras shows that this
family   depends on three parameters
$\{\bdiez\ne 0,\bcinco,\bseis\}$ with $\buno= -\bdiez$.

\subsect{Coboundary  extended Galilei bialgebras}

The next step in this procedure is to find out the extended
Galilei bialgebras that are coboundary ones. This means
that we have to deduce the classical $r$-matrices such that
\be
\delta(X)=[1\otimes X + X \otimes 1,\,  r]   \qquad  \forall
X\in \gal.
\label{bb}
\ee
It is well known that the element $r\in  \gal \otimes \gal$  
defines a coboundary Lie bialgebra $(\gal,\delta (r))$ if and
only if it fulfils the modified classical Yang--Baxter
equation (YBE)
\be
[X\otimes 1\otimes 1 + 1\otimes X\otimes 1 +
1\otimes 1\otimes X,[[r,r]]\, ]=0  \qquad \forall X\in \gal 
\label{adbc}
\ee
where $[[r,r]]$ is the Schouten bracket defined by
\be
[[r,r]]:=[r_{12},r_{13}] + [r_{12},r_{23}] + 
[r_{13},r_{23}] .
\label{adbd}
\ee
Here, if $r=r^{i j} X_i\otimes X_j$, we have denoted
$r_{12}=r^{i j} X_i\otimes X_j\otimes 1$,
$r_{13}=r^{i j} X_i\otimes 1\otimes X_j$ and
$r_{23}=r^{i j} 1\otimes X_i\otimes  X_j$. There are
 two different types of coboundary Lie bialgebras: 

\noindent
 (i) If the $r$-matrix is a skew-symmetric solution of the
classical YBE,  $[[r,r]]=0$ (the Schouten bracket vanishes),
then we obtain a {\em
non-standard} (or triangular) Lie bialgebra. 

\noindent
(ii) When  $r$
is a skew-symmetric solution of modified classical YBE
(\ref{adbc})   with non-vanishing Schouten bracket, we  
find a  {\em standard} (or quasi-triangular) Lie bialgebra.

Let us consider an arbitrary skewsymmetric element of 
$\gal \wedge \gal $:
\be
r=\aauno  K\wedge P + \aados  K\wedge M + \aatres  K\wedge H
+ \aaseis P\wedge M + \aasiete  P\wedge H + \aadiez M\wedge H.
\label{da}
\ee
The corresponding   Schouten bracket (\ref{adbd}) reads:
\bea
&&[[r,r]]= -\aatres^2  K\wedge P\wedge H +
(\aauno^2 - \aados \aatres)  K\wedge P\wedge
M \cr
&&\qquad\qquad + \aauno \aatres  K\wedge H\wedge M+
(\aauno \aasiete - \aatres \aadiez)  P\wedge H\wedge M .
\label{db}
\eea
The   modified classical  YBE (\ref{adbc}) implies
$\aatres=0$, so that the Schouten bracket reduces to
\be
[[r,r]]= \aauno^2    K\wedge P\wedge
M  + \aauno \aasiete    P\wedge H\wedge M .
\label{dc}
\ee
Hence we obtain
 a   standard classical $r$-matrix when $\aatres=0$ and
$\aauno\ne 0$, and a non-standard one when $\aatres=\aauno=0$.
 
On the other hand, the most general  element  $\eta\in \gal
\otimes \gal$ which is
$Ad^{\otimes 2}$ invariant  turns out to be
\be
\eta = \tau_1   (P\otimes P - M\otimes H - H\otimes M)
+\tau_2  M\otimes M + \tau_3  P\wedge M  
\label{dd}
\ee
where $\tau_i$ are arbitrary real numbers. Since $r'=r+ \eta$ 
generates the same bialgebra as $r$, we can choose
$\tau_1=\tau_2=0$  and $\tau_3=-\aaseis$ showing that the
term $\aaseis P\wedge M$ can be  assumed to be equal to
zero.

Both types of  coboundary  bialgebras are included in the
family Ia as follows:

\noindent
$\bullet$ {\it Standard}: $\btres= \bocho=\aauno\ne 0$, 
$\bcinco=-\aadiez$,
$\bseis=-\aasiete$, $\bdoce=-\aados$ and  $\bonce=0$.

\noindent
$\bullet$ {\it Non-standard}:  $\btres=0$, 
$\bcinco=-\aadiez$, $\bseis=-\aasiete$,
$\bdoce=-\aados$, $\bocho=0$  and  $\bonce=0$.

Furthermore the standard type  can be simplified by taking
into account the automorphism defined by
\bea
&&K'=K+\frac{\bseis}{\btres}H\qquad
H'=H-\frac{\bdoce}{\btres}P
\qquad M'=M\cr
&&P'=P-\frac{\bdoce}{\btres}M\qquad 
\bcinco'=\bcinco+\frac{\bseis\bdoce}{\btres}\qquad \btres\ne 0
\eea
which transforms the  classical $r$-matrix into
\be
r=\btres K'\wedge P' + \bcinco' H'\wedge M' +
\frac{\bcinco'\bdoce}{\btres} P'\wedge M' .
\ee
 As explained above we can discard the term $P'\wedge M'$, so
that the standard bialgebras depend on two parameters
$\{\btres\ne 0, \bcinco\}$.

For the sake of clarity the results obtained in this section
are summarized in the table 1; we display the final
cocommutators corresponding to the four families of 
bialgebras, together with   the coboundary bialgebras  as
subcases of the  family Ia.

\bigskip

{\footnotesize
\noindent
{\bf {Table 1}}.  The four multiparametric families of (1+1)
extended Galilei bialgebras.

\medskip
\noindent
\begin{tabular}{ll}
\hline
\hline
Family Ia& Six parameters:
$\{\btres,\bcinco,\bseis,\bdoce,\bocho,\bonce\}$\\[0.3cm]
&$\delta(K)=\btres K\wedge M + \bcinco P\wedge M +\bseis
H\wedge M$\\[0.1cm]
&$\delta(H)=\bonce K\wedge M + \bdoce P\wedge M
+(\bocho-\btres) H\wedge M$\\[0.1cm]
&$\delta(P)=  \bocho P\wedge M  
\qquad \delta(M)=0$\\[0.3cm]
Standard& Two parameters:
$\{\btres\ne 0,\bcinco\}$ with  $\bocho=\btres$ and
$\bseis=\bdoce=\bonce=0$\\[0.1cm]
&$r=\auno  K\wedge P+ \adiez H\wedge M $\\[0.3cm]
Non-standard& Three parameters:
$\{\bcinco,\bseis,\bdoce\}$ with 
$\btres=\bocho=\bonce=0$\\[0.1cm] 
&$r=  \adiez H\wedge M  +
 \asiete H\wedge P + \ados  M\wedge K$\\[0.3cm]
\hline
Family Ib& Three parameters:
$\{\bcuatro\ne 0,\btres,\bdoce\}$\\[0.3cm]
&$\delta(K)=\btres K\wedge M + \bcuatro P\wedge H $\\[0.1cm]
&$\delta(H)=\bdoce P\wedge M\qquad \delta(P)=\btres P\wedge
M\qquad
\delta(M)=0$\\[0.3cm]
\hline
Family IIa& One parameter:
$\{\bdiez\ne 0\}$\\[0.3cm]
&$\delta(K)=0\qquad  \delta(H)=-\bdiez P\wedge H$\\[0.1cm]
&$ \delta(P)=\bdiez
H\wedge M\qquad
\delta(M)=\bdiez P\wedge M$\\[0.3cm]
\hline
Family IIb& Three parameters:
$\{\bdiez\ne 0,\bcinco,\bseis\}$\\[0.3cm]
&$\delta(K)=-\bdiez K\wedge P + \bcinco P\wedge M
+\bseis H\wedge M$\\[0.1cm]
&$\delta(H)=0\qquad \delta(P)=0\qquad
\delta(M)=\bdiez P\wedge M$\\[0.3cm]
\hline
\hline
\end{tabular}
}


\sect{Quantum extended Galilei algebras}

We proceed to obtain the  Hopf algebras corresponding
to the four families of  (1+1) extended Galilei  bialgebras.
We shall write only the coproducts,  the compatible
commutation rules and the deformed Casimir operators; the
counit is always  trivial   and the antipode can be easily
deduced by means of the Hopf algebra axioms.


\subsect{Family Ia: quantum coboundary algebras}

All the terms appearing in the cocommutators have the form
$X\wedge M$ where $X$ is a non-primitive generator and $M$ is
primitive. Therefore we can apply the Lyakhovsky--Mudrov (LM)
formalism  \cite{LM,Lyak} in the same way as it was for the 
$h_3$, $h_4$ and $gl(2)$  algebras \cite{heis,osc,gl}
obtaining   directly the coproduct. We write the cocommutator
displayed in table 1 in matrix form as
\be
 \delta\left(\begin{array}{c}
K \\ H\\ P 
\end{array}\right)=
\left(\begin{array}{ccc}
- \btres M & -\bseis M &-\bcinco M  \\
 -\bonce M & (\btres - \bocho) M
& -\bdoce M\\
0&0& -\bocho M
\end{array}\right)\dot\wedge \left(\begin{array}{c}
K \\ H\\ P
\end{array}\right) . 
\label{ca}
\ee
Hence  the  coproduct is given by:
\be   
 \Delta\left(\begin{array}{c}
K \\ H\\ P 
\end{array}\right)=
\left(\begin{array}{ccc}
1 \otimes K  \\ 1 \otimes H \\ 1 \otimes P
\end{array}\right) +
\sigma\left( \exp\left\{\left(\begin{array}{ccc}
\btres M & \bseis M &\bcinco M  \\ \bonce M & (\bocho-\btres)
M & \bdoce M\\
0&0& \bocho M
\end{array}\right) \right\}\dot\otimes \left(\begin{array}{c}
K \\ H\\ P  
\end{array}\right)\right)  
\label{cb}
\ee
where $\sigma(X_i\otimes X_j)=X_j\otimes X_i$.
Therefore if  we denote the entries of the above matrix
exponential by $E_{ij}$,  the coproduct turns out to be:
\bea
&&\Delta(M)=1\otimes M + M\otimes 1\cr
&&\Delta(K)=1\otimes K + K\otimes E_{11}(M) + H\otimes
E_{12}(M) + P\otimes E_{13}(M)\cr
&&\Delta(H)=1\otimes H + H\otimes E_{22}(M) + K\otimes
E_{21}(M) + P\otimes E_{23}(M)\cr
&&\Delta(P)=1\otimes P + P\otimes E_{33}(M) + K\otimes
E_{31}(M) + H\otimes E_{32}(M) .
\label{cd}
\eea
The  functions $E_{ij}$ are rather complicated
and we omit them.  However as the coboundary bialgebras
belong to this family, we present in the following the
complete Hopf structure for these particular cases.

If we set $\bocho=\btres$ and
$\bseis=\bdoce=\bonce=0$  in the general expression
(\ref{cd}), we find that the   coproduct of the standard
quantum  algebra  $U^{(s)}_{\auno\ne 0,\adiez}(\gal)$ is given
by
\bea
&&\Delta(M)=1\otimes M + M\otimes 1 \qquad
\Delta(H)=1\otimes H + H\otimes 1\cr
&&\Delta(P)=1\otimes P + P\otimes e^{\auno M}\cr
&&\Delta(K)=1\otimes K + K\otimes e^{\auno M}
  +{\adiez}  P\otimes M e^{\auno M}  .
\label{ed}
\eea
The corresponding deformed commutation rules   and Casimirs
can be now obtained; they are 
\be 
 [K,H]=P
\qquad [K,P]=  \frac {e^{2\auno M}-1} {2\auno}  \qquad
[H,P]=0  \qquad
 [M,\cdot\,]=0 
\label{ee}
\ee 
\be
{\cal C}_1=M\qquad {\cal C}_2=
  P^2 - 2 \left( \frac {e^{2\btres M}-1} {2\btres}\right) H.
\label{eee}
\ee
We recall that the quantum $\gal$ algebra with
$\adiez=0$,  $U^{(s)}_{\auno\ne 0}(\gal)$,  was first
constructed in \cite{CK2d} within the framework of (1+1)
quantum Cayley--Klein algebras, and that its corresponding 
quantum deformation in (3+1) dimensions was obtained in
\cite{azcab} by means of a contraction limit of a
pseudoextension of the well-known $\kappa$-Poincar\'e algebra.

Likewise the coproduct of the non-standard quantum 
algebra  $U^{(n)}_{\adiez,\asiete,\ados}(\overline{{\cal
G}})$ comes from (\ref{cd}) provided that
$\btres=\bocho=\bonce=0$:
\bea
&&\Delta(M)=1\otimes M + M\otimes 1 \qquad
\Delta(P)=1\otimes P + P\otimes 1 \nonumber\\[2pt]
&&\Delta(H)=1\otimes H + H\otimes 1 + \ados   P\otimes M \cr
&&\Delta(K)=1\otimes K + K\otimes 1+ \adiez  P\otimes M  
+ \asiete H\otimes M 
  +\frac {\asiete \ados}{2}  P\otimes M^2 .
\label{fd}
\eea
The compatible deformed commutation rules  and Casimirs read
\be
[K,H]=P+\frac {\ados}{2}  M^2  \qquad [K,P]=  M \qquad
[H,P]=0 \qquad [M,\cdot\,]=0 
\label{fe}
\ee
\be
{\cal C}_1=M\qquad  {\cal C}_2= \left( P+ \frac {\ados}{2} 
M^2\right)^2 -2  M   H .
\label{ff}
\ee


\subsect{Family Ib: $U_{\bcuatro\ne
0,\btres,\bdoce}(\gal)$}

The presence of the term $\bcuatro P\wedge H$ in $\delta(K)$
precludes a direct use of the LM approach  since $M$ is the
only primitive generator. In spite of this fact, if  we do
not consider initially the parameter
$\bcuatro$,  the cocommutator can be written as
\be
 \delta\left(\begin{array}{c}
K \\ H\\ P 
\end{array}\right)=
\left(\begin{array}{ccc}
- \btres M & 0 &0 \\ 0 & 0
& -\bdoce M\\
0&0& -\btres M
\end{array}\right)\dot\wedge \left(\begin{array}{c}
K \\ H\\ P
\end{array}\right) . 
\label{ce}
\ee
Then the coproduct for $H$ and $P$  as well as the terms of
the coproduct of $K$ not depending on $\bcuatro$  come from
the LM method by means of:
\be   
 \Delta\left(\begin{array}{c}
K \\ H\\ P 
\end{array}\right)=
\left(\begin{array}{ccc}
1 \otimes K  \\ 1 \otimes H \\ 1 \otimes P
\end{array}\right) +
\sigma \left( \exp\left\{\left(\begin{array}{ccc}
 \btres M & 0 &0 \\ 0 & 0
&  \bdoce M\\
0&0&  \btres M
\end{array}\right)\right\}\dot\otimes \left(\begin{array}{c}
K \\ H\\ P  
\end{array}\right)\right) .
\label{cf}
\ee
The remaining terms of the coproduct of $K$ (whose first
order in the bialgebra parameters lead to  
$\bcuatro P\wedge H$ in $\delta(K)$) can be computed 
by solving the coassociativity condition. The resultant  
coproduct for the three-parameter quantum 
 algebra $U_{\bcuatro\ne
0,\btres,\bdoce}(\gal)$
reads  
\bea
&&\Delta(M)=1\otimes M + M \otimes 1\qquad
\Delta(P)=1\otimes P + P \otimes e^{\btres M}\nonumber\\[2pt]
&&\Delta(H)=1\otimes H + H \otimes 1 + \bdoce P\otimes
\left(\frac {e^{\btres M}-1}{\btres}\right)\label{cg}\\
&&\Delta(K)=1\otimes K + K \otimes e^{\btres M}   + \bcuatro
P\otimes H e^{\btres M} +\frac {\bcuatro \bdoce}2 P^2\otimes 
\left(\frac {e^{\btres M}-1}{\btres}\right) e^{\btres M}.
\nonumber
\eea
The compatible commutation rules can be now deduced
\bea
&&[K,H]=P+ \frac{\bdoce}{2}\left(\frac {e^{\btres
M}-1}{\btres} \right)^2
\qquad [H,P]=0 \cr
&& [K,P]=\frac {e^{2\btres M}-1} {2\btres}  \qquad \quad
[M,\cdot\,]=0,
\label{ch}
\eea
and the quantum Casimirs read
\be
{\cal C}_1=M\qquad {\cal C}_2=
\left( P+ \frac{\bdoce}{2}\left(\frac {e^{\btres
M}-1}{\btres} \right)^2\right)^2 - 2 \left( \frac {e^{2\btres
M}-1} {2\btres}\right) H.
\label{ci}
\ee


\subsect{Family IIa}

We denote by $\{k,h,p,m\}$  the generators spanning the dual
basis of $\gal$. The  one-parameter  bialgebra written in
table 1 allows us to obtain  the  following dual Lie brackets 
\be
[p,m]=\bdiez m\qquad [p,h]=-\bdiez h\qquad
[h,m]=\bdiez p \qquad [k,\cdot\,]=0  
\label{ccx}
\ee 
which close a  Lie algebra isomorphic to $gl(2)$; $k$ plays
the role of the central generator. Hence the group law of
$GL(2)$ arises as the coproduct of the quantum  algebra of
this  family IIa. This fact is in agreement with the
classification of $gl(2)$ bialgebras carried out in
\cite{gl}; conversely, it can be checked that the coproduct
of the non-standard family II of quantum $gl(2)$ algebras
with   $b_-=b=0$ constructed in \cite{gl}, 
$U_{b_+}(gl(2))$,  is the group law of the (1+1) extended
Galilei group.


\subsect{Family IIb: $U_{\bdiez\ne 0,\bcinco,\bseis}(\gal)$}

In this case $H$ and $P$ are two commuting primitive
generators, so that  the cocommutators of the two remaining
generators can be expressed as
\be
 \delta\left(\begin{array}{c}
K \\ M 
\end{array}\right)=
\left(\begin{array}{cc}
0 &\bseis H   \\ 0 &0
\end{array}\right)\dot\wedge \left(\begin{array}{c}
K \\ M 
\end{array}\right)+
\left(\begin{array}{cc}
\bdiez P & \bcinco P   \\ 0 & \bdiez P 
\end{array}\right)\dot\wedge \left(\begin{array}{c}
K \\ M
\end{array}\right)  .
\label{cp}
\ee
Therefore their coproduct is  provided by the LM method:
\be   
 \Delta\left(\begin{array}{c}
K \\ M
\end{array}\right)=
\left(\begin{array}{ccc}
1 \otimes K  \\ 1 \otimes M
\end{array}\right) +
\sigma\left( \exp\left\{\left(\begin{array}{cc}
 -\bdiez P & -\bcinco P - \bseis H      \\ 0
  &  -\bdiez P
\end{array}\right) \right\}\dot\otimes \left(\begin{array}{c}
K \\ M 
\end{array}\right)\right) .
\label{cq}
\ee
Hence the explicit coproduct of the three-parameter quantum 
algebra $U_{\bdiez\ne 0,\bcinco,\bseis}(\gal)$ turns out to be
\bea
&&\Delta(P)=1\otimes P + P\otimes 1\qquad 
\Delta(H)=1\otimes H + H\otimes 1 \cr
&&\Delta(M)=1\otimes M + M\otimes e^{-\bdiez P}\cr
&&\Delta(K)=1\otimes K + K\otimes e^{-\bdiez P}
-M\otimes (\bcinco P + \bseis H)  e^{-\bdiez P} .
\label{cr}
\eea
The compatible deformed commutation rules are given by:
\bea
&&[K,H]=\frac{1- e^{-\bdiez P}}{\bdiez} \qquad
[K,P]=M\qquad [H,P]=0\cr
&&[M,K]= \frac{1}{2}\bdiez M^2\qquad [M,H]=0\qquad [M,P]=0 ,
\label{ct}
\eea
while the deformed Casimir operators read
\be
{{\cal C}_1}=e^{\bdiez P/2} M
\qquad
{{\cal C}_2}=\left( \frac{\sinh(\bdiez
P/4)}{\bdiez/4}\right)^2
 -2 e^{\bdiez P/2} M H.
\label{cu}
\ee

The particular quantum deformation with
$\bcinco=\bseis=0$,   $U_{\bdiez\ne 0}(\gal)$, was originally
obtained in \cite{ita,itab}. More explicitly, it can be
checked that  the generators $\{{\cal B},{\cal T},{\cal
P},{\cal M}\}$ and deformation parameter $a$ defined by
\be
{\cal B}=i e^{\bdiez P/2} K\qquad
{\cal T}=iH\qquad
{\cal P}=iP\qquad
{\cal M}=i e^{\bdiez P/2} M\qquad
a=\bdiez /2
\ee
give rise to the quantum extended Galilei algebra introduced
in \cite{ita,itab}. We recall that this quantum algebra
$U_{a\ne 0}(\gal)$ was shown to describe the symmetry of
magnons on the one-dimensional Heisenberg ferromagnet for both
the isotropic (XXX) and the anisotropic (XXZ) magnetic
chain; the quantum algebra symmetry was completely equivalent
to the Bethe ansatz  and the deformation parameter was
identified with the chain spacing.


\sect{Quantum universal $R$-matrices}

In this section we deduce  quantum universal $R$-matrices
associated to the standard and non-standard quantum extended
Galilei algebras obtained within  the family Ia in the
section 3.1.

\subsect{Standard universal $R$-matrix}

We consider the standard classical $r$-matrix
\be
r=\auno  K\wedge P+ \adiez H\wedge M\qquad  \auno\ne 0.
\label{ma}
\ee
If we look for  a non-skewsymmetric classical $r$-matrix by
adding a generic $Ad^{\otimes 2}$ invariant  element
$\eta$ (\ref{dd}) to (\ref{ma}) and we impose the
classical YBE   to be fulfilled   we find that the parameter
$\auno$ must be equal to zero. Consequently there does not
exist a quasitriangular universal $R$-matrix  satisfying the
quantum YBE, whose first order in the deformation parameters
gives  the standard $r$-matrix (\ref{ma}). However, as we
shall show in the following, it is possible to find a
non-quasitriangular universal $R$-matrix once we set
$\adiez=0$.

The coproduct  and commutation rules  of the
one-parameter quantum algebra  $U^{(s)}_{\auno\ne 0}(\gal)$
are obtained from (\ref{ed}) and  (\ref{ee}) provided that
$\adiez=0$:
\bea
&&\Delta(M)=1\otimes M + M\otimes 1 \qquad 
\Delta(H)=1\otimes H + H\otimes 1 \cr
&&\Delta(P)=1\otimes P + P\otimes e^{\auno M} \qquad
 \Delta(K)=1\otimes K + K\otimes e^{\auno M}  ,
\label{mb}
\eea
\be 
 [K,H]=P  \qquad [K,P]=\frac {e^{2\auno M}-1}{2\auno} \qquad
[H,P]=0 \qquad
 [M,\cdot\,]=0.
\label{mc}
\ee 
The crucial point  now is that the three generators $K$, $P$
and $M$ close a  Hopf subalgebra deforming a
Heisenberg  algebra which can be easily related with the
non-quasitriangular quantization of the Heisenberg algebra
developed in \cite{RHeisenberg} by means of
\be
K\to K e^{-\auno M/2}\qquad  
P\to P e^{-\auno M/2}\qquad 
M\to M . 
\ee
Therefore the universal  $R$-matrix given in 
\cite{RHeisenberg} which is not a solution of the quantum YBE
but it verifies
\be
{\cal R}\Delta(X){\cal R}^{-1}=\sigma\circ\Delta(X),
\label{md}
\ee
can be adapted to our basis as
\bea
&&{\cal R}=\exp(\auno K\wedge P f(M,\auno)) \nonumber\\[2pt]
&&f(M,\auno)=\frac{e^{-\auno M/2}\otimes e^{-\auno
M/2}}{\sqrt{\sinh \auno M\otimes \sinh \auno
M}}\arcsin\left({\frac{\sqrt{\sinh \auno M\otimes
\sinh \auno M}}{\cosh((\auno/2)\Delta(M))}}\right) .
\label{me}
\eea
Furthermore, it is straightforward to prove that the
relation (\ref{md}) also holds for the remaining generator
$H$, thus we conclude that (\ref{me}) is a
non-quasitriangular  universal $R$-matrix for
$U^{(s)}_{\auno\ne 0}(\gal)$.


\subsect{Non-standard universal $R$-matrix}

The non-standard classical $r$-matrix is given by
\be
r=  \adiez H\wedge M  +
 \asiete H\wedge P + \ados  M\wedge K .
\label{mf}
\ee
The corresponding universal $R$-matrix which  satisfies the
property (\ref{md}) for the {\it whole} family
$U^{(n)}_{\adiez,\asiete,\ados}(\overline{{\cal G}})$, with
coproduct (\ref{fd}) and commutation rules (\ref{fe}), turns
out to be
\be
 {\cal R} =\exp\{ - \adiez  M\otimes H + \ados  M\otimes K\}
\exp\{\asiete  H\wedge P\}
\exp\{ \adiez  H\otimes M-\ados  K\otimes M \}.
\label{mg}
\ee
The proof for $M$ is trivial since it is a primitive and
central generator.  We summarize the main steps of the
computations for the  remaining generators. If we denote
(\ref{mg}) as ${\cal R}=e^{A_3}e^{A_2}e^{A_1}$, then we find
that 
\bea
&&e^{A_1}\Delta(P)e^{-A_1}=
1\otimes P + P\otimes 1 - \ados M\otimes M  \equiv h \cr
&&e^{A_2}he^{-A_2}= h \qquad
 e^{A_3}he^{-A_3}= \sigma\circ\Delta(P).
\label{mh}
\eea
\bea
&&e^{A_1}\Delta(H)e^{-A_1}=
1\otimes H + H\otimes 1 -\frac{\ados^2}{2}(M^2\otimes M +
M\otimes M^2)\equiv f \cr
&&e^{A_2}fe^{-A_2}= f \qquad
 e^{A_3}fe^{-A_3}= \sigma\circ\Delta(H).
\label{mi}
\eea
\bea
&&e^{A_1}\Delta(K)e^{-A_1}= 1\otimes K + K\otimes 1 + \asiete
H\otimes M -\frac{\asiete\ados}2  P\otimes M^2\cr
&&\qquad\quad - \frac{\adiez\ados}{2}(M^2\otimes M + M\otimes
M^2)- \frac{\asiete\ados^2}{2} M^2\otimes M^2\equiv g_1 \cr
&&e^{A_2}g_1e^{-A_2}=
1\otimes K + K\otimes 1 + \asiete M\otimes H
-\frac{\asiete\ados}2 
 M^2\otimes P\cr
&&\qquad\quad - \frac{\adiez\ados}{2}(M^2\otimes M + M\otimes
M^2)- \frac{\asiete\ados^2}{2} M^2\otimes M^2\equiv g_2 \cr
&&e^{A_3}g_2e^{-A_3}= \sigma\circ\Delta(K).
\label{mj}
\eea
 
The question of whether (\ref{mg}) is a solution of
the quantum YBE remains as an open problem.


\sect{Classical integrable systems from Poisson $\gal$
coalgebras}

We regard now the commutation rules (\ref{bc}) as Poisson
brackets and we consider the
usual one-particle  phase space representation $D$ of $\gal$
given by
\bea
&&f_P^{(1)}=D(P)=p_1\qquad
f_M^{(1)}=D(M)=m_1\cr
&&f_K^{(1)}=D(K)=m_1 q_1\qquad
f_H^{(1)}=D(H)=\frac{p_1^2}{2 m_1}
\label{oa}
\eea
where $m_1$ is a real constant. The realization of the
Casimirs (\ref{bd}) is $C_1^{(1)}=D({\cal C}_1)=m_1$ and 
 $C_2^{(1)}=D({\cal C}_2)=0$.

The $\gal$ Lie--Poisson algebra is endowed with a Poisson
coalgebra structure by means of the primitive coproduct:
$\Delta(X)=1\otimes X+ X\otimes 1$;
this leads to  
two-particle phase space functions obtained as
$f_X^{(2)}=(D\otimes D)(\Delta(X))$:
\bea
&&f_P^{(2)}=p_1+ p_2\qquad
f_M^{(2)}=m_1+m_2 \cr
&&f_K^{(2)}=m_1 q_1+m_2 q_2\qquad
f_H^{(2)}=\frac{p_1^2}{2 m_1}+\frac{p_2^2}{2 m_2},
\label{ob}
\eea
which close again a $\gal$ algebra 
with respect to the usual Poisson bracket
$\{q_i,p_j\}=\delta_{ij}$. The formalism developed in
\cite{BR} ensures that the two-particle Hamiltonian
${\cal H}^{(2)}$ defined as the coproduct of any smooth
function ${\cal H}(K,H,P,M)$ of the coalgebra generators
\be
{\cal H}^{(2)}=(D\otimes D)(\Delta({\cal H}))={\cal
H}\bigl(f_{K}^{(2)},f_{H}^{(2)},f_{P}^{(2)},f_{M}^{(2)}\bigr) 
\label{oc}
\ee
is completely integrable. Its integral of motion
is provided by the $D\otimes D$ representation of the
coproduct of the second-order Casimir and reads:
\be
C^{(2)}_2=(D\otimes D)(\Delta({\cal C}_2))=-\frac{(m_2 p_1 -
m_1 p_2)^2}{m_1 m_2} ,
\label{od}
\ee
while the Casimir ${\cal
C}_1=M$ gives rise to a trivial integral of motion: 
$C^{(2)}_1=m_1+m_2$. A particular subset of integrable
Hamiltonians can be found by setting 
\be
{\cal H}=H+  {\cal F}(K)
\label{hamil}
\ee
 where ${\cal F}$ is any smooth function
of the boost $K$. In this case, (\ref{oc}) leads to the
natural two-particle Hamiltonian
\be
{\cal H}^{(2)}=f_{H}^{(2)}+  {\cal
F}\bigl(f_{K}^{(2)}  \bigr) =
\frac{p_1^2}{2 m_1}+\frac{p_2^2}{2 m_2}+
{\cal F}\bigl(m_1 q_1+m_2 q_2\bigr),
\label{oe}
\ee
that, by construction,  Poisson-commutes with (\ref{od}).
We stress that the generalization to integrable $N$-particle
systems can be obtained by making use of higher order
coproducts \cite{BR}.

The very same procedure can be carried out with the quantum
$\gal$ coalgebras obtained in section 3: if  we consider
the deformed commutation rules   as   Poisson brackets then
the coproduct defines the (deformed) coalgebra structure.
Therefore, once   a  one-particle
 phase space representation is deduced for each Poisson
deformed $\gal$ coalgebra,  the coproduct defines the 
two-particle phase space functions which automatically fulfil
the corresponding (deformed) Poisson brackets. In this way,
any function of the deformed generators gives rise to a
 completely integrable Hamiltonian whose  integral of motion is
again  given by the coproduct of the deformed Casimir.
All the information needed to construct  two-particle
integrable systems is displayed in table 2;  for each
multiparametric Poisson 
$\gal$ coalgebra 
 we  write its corresponding one- and two-particle phase
realization ($f_X^{(1)}$ and $f_X^{(2)}$) together with the
integrals of motion  $C^{(2)}_1$ and $C^{(2)}_2$; for all of
them, the one-particle Casimirs are $C^{(1)}_1=m_1$ and
$C^{(1)}_2=0$. We have also introduced a `deformed mass
function' defined by
\be
\masa_i(x):=\frac{e^{x m_i}-1}{x}\qquad i=1,2
\label{masa}
\ee
where $x$ is a deformation parameter  (either $\auno$ or
$2\auno$); obviously, $\lim_{x\to 0}\masa_i(x)=m_i$.

In this context, the different quantum deformations of $\gal$
can be interpreted as the structures generating
multiparametric integrable deformations of the Hamiltonians
coming from
${\cal H}$ functions. For instance, let us consider again the
Hamiltonian (\ref{hamil}) with $H$ and $K$ being
now the (Poisson) generators of deformed Galilei coalgebras.
 When ${\cal H}$ is defined on the standard
Poisson coalgebra
$U^{(s)}_{\auno\ne 0,\adiez}(\gal)$, the Hamiltonian
(\ref{oe}) is deformed into (see standard family Ia  in table
2)
\bea
&&\!\!\!\!\!\!{\cal H}_{\auno\ne 0,\adiez}^{(2)}=
f_{H}^{(2)}+  {\cal
F}\bigl(f_{K}^{(2)}  \bigr) \cr
&&\quad =
\frac{p_1^2}{2 \masa_1(2\auno)} +
\frac{p_2^2}{2 \masa_2(2\auno)}+
{\cal F}\bigl(e^{\auno m_2}\masa_1(2\auno) q_1 +
\masa_2(2\auno) q_2+\adiez e^{\auno m_2} m_2 p_1\bigr) 
\label{1ast}
\eea
which is in involution with the corresponding coproduct
of the deformed Casimir, namely
\be
C^{(2)}_2=-\frac{\left(\masa_2(2\auno) p_1 -
\masa_1(2\auno) e^{\auno m_2}p_2\right)^2}{\masa_1(2\auno)
\masa_2(2\auno)}.
\label{1astc}
\ee
Notice that the deformation parameter $\adiez$ induces a
$p_1$-dependent term in the potential. If $\adiez=0$, we see
that (\ref{1ast}) is an integrable deformation of (\ref{oe})
in which both masses have been deformed, $m_i\rightarrow
\masa_i$, and the potential is an arbitrary
function of  
$(\alpha_1\,\masa_1\, q_1 +
\masa_2\,q_2)$, where the constant $\alpha_1$ has to be
exactly
$e^{\auno m_2}$. This result can be extended to arbitrary
dimension by following \cite{BR}  (see also
\cite{chains} for the construction of integrable systems
associated to non-standard  Poisson
 $sl(2,\R)$  coalgebras).
That procedure leads to a Hamiltonian
of the type
\be
{\cal H}_{\auno\ne 0,\adiez=0}^{(N)}=
\sum_{i=1}^{N}{\frac{p_i^2}{2 \masa_i}} +
{\cal F}\bigl(\alpha_1\,\masa_1\,q_1 + \alpha_2\,\masa_2\,q_2
+\dots+\masa_N\,q_N\bigr) 
\label{1astN}
\ee
where the deformed masses and constants are
\be
\masa_i=\masa_i(2\auno)\quad i=1,\dots,N 
\qquad \alpha_l=e^{\auno (m_{l+1}+m_{l+2}+\dots+m_{N})} \quad
 l=1,\dots,N-1 .
\label{const}
\ee
The $(N-1)$ integrals
of the motion in involution with (\ref{1astN}) would be
obtained through the $k$-th coproducts
$(k=2,\dots,N)$ of the Casimir ${\cal C}_2$. 

From table 2, it is easy to check that integrable deformations
generated by the non-standard Poisson coalgebra
$U^{(n)}_{\adiez,\asiete,\ados}(\overline{{\cal G}})$ provide
only additional terms depending on $p_1$ with respect to the
non-deformed construction. Next, the family Ib 
$U_{\bcuatro\ne 0,\btres,\bdoce}(\gal)$ encompasses
simultaneously  properties of the two previous families.
Finally, the family IIb $U_{\bdiez\ne 0,\bcinco,\bseis}(\gal)$
gives rise to an essentially different integrable deformation;
if we consider again the same dynamical Hamiltonian ${\cal H}$
(\ref{hamil}) we find (for the particular case with
$\bcinco=\bseis=0$)
\bea
&&\!\!\!\!\!\!{\cal H}_{\bdiez\ne 0,\bcinco=\bseis=0}^{(2)}=
f_{H}^{(2)}+  {\cal
F}\bigl(f_{K}^{(2)}  \bigr) \cr
&&\qquad\qquad\quad=
\frac {1}{2 m_1} 
\left(\frac{\sinh(\bdiez p_1/4)}{\bdiez/4}
\right)^2+\frac {1}{2 m_2} 
\left(\frac{\sinh(\bdiez p_2/4)}{\bdiez/4}
\right)^2 \cr
&&\qquad\qquad\qquad + 
{\cal F}\left( m_1 e^{-\bdiez p_1/2 }e^{-\bdiez p_2} q_1+
m_2 e^{-\bdiez p_2/2 }q_2 \right).
\label{iib}
\eea
Hence a deformation of the kinetic energy in terms of
hyperbolic functions 
is obtained, and the
potential  is also deformed through exponentials of
the momenta. As expected, the hyperbolic functions of $p_i$
are also present in the deformed integral of the motion (see
$C_2^{(2)}$ in table 2).

\newpage

\bigskip

{\footnotesize
\noindent
{\bf {Table 2}}. Two-particle integrable systems from Poisson
$\gal$ coalgebras.

\medskip
\noindent
\begin{tabular}{l}
\hline
\hline
Family Ia: Standard Poisson coalgebra $U^{(s)}_{\auno\ne
0,\adiez}(\gal)$\\[0.2cm] 
$\displaystyle{f_K^{(1)}=\masa_1(2\auno) q_1\qquad
f_H^{(1)}=\frac{p_1^2}{2 \masa_1(2\auno)}\qquad
f_P^{(1)}=p_1\qquad
f_M^{(1)}=m_1}$\\[0.3cm]
$\displaystyle{f_K^{(2)}=e^{\auno m_2}\masa_1(2\auno) q_1 +
\masa_2(2\auno) q_2+\adiez e^{\auno m_2} m_2 p_1
\qquad    f_M^{(2)}=m_1+m_2}$\\[0.2cm]
$\displaystyle{
f_H^{(2)}=\frac{p_1^2}{2 \masa_1(2\auno)} +
\frac{p_2^2}{2 \masa_2(2\auno)}\qquad
f_P^{(2)}= e^{\auno m_2} p_1+
p_2}$\\[0.3cm]
$\displaystyle{
C_1^{(2)}=m_1+m_2\qquad C^{(2)}_2=-\frac{\left(\masa_2(2\auno) p_1 -
\masa_1(2\auno) e^{\auno m_2}p_2\right)^2}{\masa_1(2\auno)
\masa_2(2\auno)} }$\\[0.3cm]
\hline
Family Ia: Non-standard Poisson coalgebra
$U^{(n)}_{\adiez,\asiete,\ados}(\overline{{\cal G}})$\\[0.2cm] 
$\displaystyle{f_K^{(1)}=m_1 q_1\qquad
f_H^{(1)}=\frac{p_1^2}{2 m_1}\qquad
f_P^{(1)}=p_1 - \frac{\ados}{2}m_1^2\qquad
f_M^{(1)}=m_1}$\\[0.3cm]
$\displaystyle{f_K^{(2)}=m_1  q_1 +
m_2  q_2+ \left(p_1 - \frac{\ados}{2} m_1^2\right
) \left(\adiez m_2 + \frac 12 \asiete\ados m_2^2 \right) 
+\asiete m_2 \frac{p_1^2}{2 m_1} \qquad f_M^{(2)}=m_1+m_2}$\\[0.3cm]
$\displaystyle{f_H^{(2)}=\frac{p_1^2}{2 m_1} +
\frac{p_2^2}{2 m_2} + \ados m_2\left(p_1 - \frac{\ados}{2} m_1^2\right
) \qquad f_P^{(2)}=  p_1+
p_2 - \frac{\ados}{2}(m_1^2+m_2^2) }$\\[0.3cm]
$\displaystyle{C_1^{(2)}=m_1+m_2\qquad C^{(2)}_2=-\frac{(m_2 p_1 - m_1
p_2)^2}{m_1 m_2} + 2 \ados m_2 (m_1 p_2 - m_2 p_1)+ \ados^2 m_1^2 m_2
(m_1 + 2 m_2) }$\\[0.3cm]
\hline
Family Ib: $U_{\bcuatro\ne
0,\btres,\bdoce}(\gal)$\\[0.2cm] 
$\displaystyle{f_K^{(1)}=\masa_1(2 \btres) q_1\qquad
f_H^{(1)}=\frac{p_1^2}{2 \masa_1(2 \btres) }\qquad
f_P^{(1)}=p_1 - \frac{\bdoce}{2}\masa_1^2(\btres) \qquad
f_M^{(1)}=m_1}$\\[0.2cm]
$\displaystyle{f_K^{(2)}=e^{\btres m_2}\masa_1(2 \btres) q_1+
\masa_2(2 \btres) q_2+\bcuatro e^{\btres m_2}\left(p_1 
- \frac{\bdoce}{2}\masa_1^2(\btres)\right)
\frac{p_2^2}{2 \masa_2(2 \btres) } }$\\[0.2cm]
$\displaystyle{\qquad\qquad  +
\frac{\bcuatro \bdoce}{2}e^{\btres m_2}\left( p_1 
- \frac{\bdoce}{2}\masa_1^2(\btres)\right)^2  \masa_2(\btres)
 }$\\[0.2cm]
$\displaystyle{f_H^{(2)}=
\frac{p_1^2}{2 \masa_1(2 \btres) }+\frac{p_2^2}{2 \masa_2(2 \btres) }
+\bdoce\left(p_1 
- \frac{\bdoce}{2}\masa_1^2(\btres)\right)\masa_2(\btres)\qquad
 f_M^{(2)}=m_1+m_2 }$\\[0.2cm]
$\displaystyle{f_P^{(2)}=e^{\btres m_2} p_1 +p_2 -\frac{\bdoce}{2}
\left( \masa_1^2(\btres)e^{\btres m_2} + \masa_2^2(\btres)\right)
\qquad C_1^{(2)}=m_1+m_2 }$\\[0.2cm]
$\displaystyle{C^{(2)}_2= -\frac{\left(\masa_2(2\auno) p_1 -
\masa_1(2\auno) e^{\auno m_2}p_2\right)^2}{\masa_1(2\auno)
\masa_2(2\auno)} - 2  \bdoce \masa_2(\auno)\left(\masa_2(2\auno)p_1
-\masa_1(2\auno) e^{\auno m_2}p_2
\right)}$\\[0.2cm]
$\displaystyle{
\qquad\quad +\frac{\bdoce^2}{4}\masa_1^2(\auno)\masa_2(\auno)\bigl\{ 
\masa_2(\auno)\left(2+\btres \masa_1(\auno)\right)^2
\left(1+\btres \masa_2(\auno)\right)^2 
}$\\[0.2cm]
$\displaystyle{\qquad\qquad\qquad\qquad \qquad\qquad 
  + 4\masa_1(2\auno)
+4 \masa_2(2\auno)+8\auno \masa_1(2\auno)\masa_2(2\auno) \bigr\}
}$\\[0.2cm]
\hline
Family IIb: $U_{\bdiez\ne 0,\bcinco,\bseis}(\gal)$\\[0.2cm] 
$\displaystyle{f_K^{(1)}=m_1 e^{-\bdiez p_1/2 }q_1\qquad
f_H^{(1)}=\frac 1{2 m_1}\left(\frac{\sinh(\bdiez p_1/4)}{\bdiez/4}
\right)^2\qquad f_P^{(1)}=p_1 \qquad
f_M^{(1)}= e^{-\bdiez p_1/2 } m_1}$\\[0.2cm]
$\displaystyle{
f_K^{(2)}= m_1 e^{-\bdiez p_1/2 }e^{-\bdiez p_2} q_1+
m_2 e^{-\bdiez p_2/2 }q_2 -m_1  e^{-\bdiez p_1/2 }e^{-\bdiez
p_2}\left(
\bcinco p_2 + \frac {\bseis }{2 m_2} 
\left(\frac{\sinh(\bdiez p_2/4)}{\bdiez/4}
\right)^2 \right)
}$\\[0.3cm]
$\displaystyle{f_H^{(2)}=\frac {1}{2 m_1} 
\left(\frac{\sinh(\bdiez p_1/4)}{\bdiez/4}
\right)^2+\frac {1}{2 m_2} 
\left(\frac{\sinh(\bdiez p_2/4)}{\bdiez/4}
\right)^2 \qquad f_P^{(2)}=p_1+p_2
}$\\[0.3cm] 
$\displaystyle{ f_M^{(2)}=m_1 e^{-\bdiez p_1/2 }e^{-\bdiez p_2}+
m_2 e^{-\bdiez p_2/2 }\qquad
C_1^{(2)}=m_1 e^{-\bdiez p_2/2} +m_2 e^{\bdiez p_1/2}
}$\\[0.2cm] 
$\displaystyle{ 
C_2^{(2)}=-\frac {1}{m_1 m_2} \left( m_2
\left(\frac{\sinh(\bdiez p_1/4)}{\bdiez/4}
\right) e^{\bdiez p_1/4}- m_1  \left(\frac{\sinh(\bdiez
p_2/4)}{\bdiez/4}
\right) e^{-\bdiez p_2/4}
\right)^2}$\\[0.3cm] 
\hline
\hline
\end{tabular}
}

\bigskip


\sect{Concluding remarks}

To end with, we would like to comment on the relationship
between $\gal$ and the (1+1)-dimensional free
heat-Schr\"odinger equation (HSE). This can be established by
recalling the usual kinematical differential realization of
the  Galilei generators in terms of the space and time 
coordinates
$(x,t)$:
\be
K=-t\partial_x - m  x \qquad H=\partial_t\qquad
P=\partial_x\qquad M=m
\label{na}
\ee
where  $m$ (the mass) is a constant that labels the
representation. The action of the Casimir ${\cal C}_2$
(\ref{bd}) on a function $\Psi(x,t)$ through (\ref{na}) gives
rise to the  HSE:
\be
\left\{\partial_x^2 -2 m
\partial_t \right\}\Psi(x,t)=0 .
\label{nb}
\ee
The quantum $\gal$ algebras obtained in section 3 allow us to
deduce in a straightforward way  deformed HSE's by following
a similar procedure to the non-deformed case. In particular,
once  a deformed differential representation is found for 
each multiparametric quantum  $\gal$ algebra, the  deformed
HSE is provided by  the quantum Casimir written in terms of
such a  representation; hence  the resulting HSE has
automatically a quantum $\gal$ algebra symmetry.
In particular, if we consider the quantum algebra 
$U_{\bdiez\ne 0,\bcinco,\bseis}(\gal)$ of the family IIb, we
find the following differential-difference realization:  
\be
K=-t\left(\frac{1-e^{- \bdiez \partial_x}}{\bdiez}\right)
-m x e^{-\bdiez \partial_x/2} \qquad H=\partial_t\qquad
P=\partial_x\qquad M=m  e^{-\bdiez \partial_x/2} .
\label{wa}
\ee
Hence we obtain a  space discretized HSE in a uniform
lattice with $U_{\bdiez\ne 0,\bcinco,\bseis}(\gal)$ symmetry
given by
\be
\left\{\left( 
\frac{\sinh(\bdiez \partial_x/4)}{\bdiez/4}
\right)^2 - 2 m \partial_t 
\right\}\Psi(x,t)=0 .
\label{wb}
\ee
It can be easily checked that  the remaining quantum $\gal$
algebras would also lead to `deformed' equations but with no
discretization. Finally, we recall that a similar equation to
(\ref{wb}) with quantum Schr\"odinger algebra symmetry has
been obtained in \cite{equation}.


{\section*{Acknowledgments}}
\noindent
 A.B. and
F.J.H. have been partially supported by  Junta
de Castilla y Le\'on, Spain  (Project   CO2/399).

\end{document}